\documentclass{svproc}

\usepackage{url}

\usepackage{fixltx2e}

\usepackage[ansinew]{inputenc}
\usepackage{srcltx}

\usepackage{datetime}
\usepackage{ngerman}
\usepackage{textcomp}

\usepackage{graphicx, epsfig}
\usepackage{tikz}
\usepackage{tkz-graph}

\usepackage{amsmath}

\usepackage{amssymb}
\usepackage{stmaryrd}
\usepackage{bm}

\usepackage{enumerate}

\usepackage{afterpage}
\usepackage{xspace}
\usepackage{listings}

\lstdefinelanguage{Sage}[]{Python}
{morekeywords={True,False,sage,singular},
sensitive=true}
\definecolor{dblackcolor}{rgb}{0.0,0.0,0.0}
\definecolor{dbluecolor}{rgb}{.01,.02,0.7}
\definecolor{dredcolor}{rgb}{0.8,0,0}
\definecolor{dgraycolor}{rgb}{0.30,0.3,0.30}

\newcommand\´{\kern 1pt}

\newcommand\noms{\hspace{-\mathsurround}}

\newcommand\Rand[1]{
  \marginpar{\raggedleft\scriptsize\hspace{0pt}#1}}%

\renewcommand{\(}{\noms$}
\renewcommand{\)}{\noms$}
\renewcommand\frac[2]{\genfrac{}{}{.4pt}{}{#1}{#2}}

\renewcommand\dfrac[2]{\genfrac{}{}{.4pt}{0}{#1}{#2}}

\newcommand\mathRand[1]{\hspace{\mathsurround}\Rand{#1}\nolinebreak\noms}
\def\rand #1"#2"{\mathRand{\(#2\)}#1#2}
\def\randd #1"#2"#3\randd#4"#5"{\mathRand{\(#2\), \(#5\)}#1#2#3#4#5}

\newcommand\eqby[2][=]%
  {\ensuremath{\overset{\makebox[0pt]{\ensuremath{\smash[t]{\scriptstyle#2}}}}{#1}}}

\newcommand\F[1][\ ]{\mathbb{F}_{\!#1}}

\newcommand\Rl{\mathbb{R}}
\newcommand\Q{\mathbb{Q}}

\newcommand\N{\mathbb{N}}

\newcommand\nin{\notin}

\newcommand\ssm{\textup{\texttt{\char92}}}

\newcommand\mto{\mapsto}

\newcommand\lto{\longrightarrow}

\newcommand\lTo{\Longrightarrow}

\newcommand\nach{\mathbin\circ}

\newcommand\ex{\exists\,}

\newcommand\DP{\colon\discretionary{\!\kern -.17em}{}{}}
\newcommand\mitsymbol{\textup{\textbrokenbar}}
\renewcommand\mit{\,\ \discretionary{\mitsymbol}{}{}\mitsymbol\ \,}

\renewcommand\div{\mathrel{\bigm\lfloor\!\!\!\bigm\lfloor}}
\newcommand\vid{\mathrel{\bigm\rfloor\!\!\!\bigm\rfloor}}
\newcommand\ndiv{\mathrel{\;\!\div\hspace{-12pt}\kern0pt\lower2pt%
  \hbox{\ensuremath{^\diagup}}\!}}
\newcommand\ndivps{\mathrel{\;\!\div\hspace{-9pt}\kern0pt\lower2pt%
  \hbox{\ensuremath{^\diagup}}\!}}

\newcommand\nvid{\mathrel{\;\!\vid\hspace{-12pt}\kern0pt\lower2pt%
  \hbox{\ensuremath{^\diagup}}\!}}
\newcommand\nvidps{\mathrel{\;\!\vid\hspace{-9pt}\kern0pt\lower2pt%
  \hbox{\ensuremath{^\diagup}}\!}}
\providecommand\abs[1]{\lvert#1\rvert}
\providecommand\Abs[1]{\bigl\lvert#1\bigr\rvert}

\newcommand\LECC{List Edge Coloring Conjecture\xspace}
\newcommand\CN{Combinatorial Nullstellensatz\xspace}

\newcommand\ä{\alpha}
\renewcommand\d{\delta}
\DeclareMathOperator\sgn{sgn}
\DeclareMathOperator\intt{int}

\DeclareMathOperator\OF{OF}

\newcommand\cop[1]{\lfloor#1\rfloor}

\begin{document}
\mainmatter              
\title{Orientations of 1-Factorizations and\\ the List Chromatic Index of Small Graphs}
\titlerunning{Orientations of 1-Factorizations}  
%
\author{Uwe Schauz}
\authorrunning{Uwe Schauz} 
%
\tocauthor{Uwe Schauz}
\institute{Xi’an Jiaotong-Liverpool University, Suzhou 215123, China,\\
\email{uwe.schauz@xjtlu.edu.cn}
}

\maketitle              

\begin{abstract}
As starting point, we formulate a corollary to the Quantitative Combinatorial Nullstellensatz. This
corollary does not require the consideration of any coefficients of polynomials, only evaluations
of polynomial functions. In certain situations, our corollary is more directly applicable and more
ready-to-go than the Combinatorial Nullstellensatz itself. It is also of interest from a numerical
point of view. We use it to explain a well-known connection between the sign of 1-factorizations
(edge colorings) and the List Edge Coloring Conjecture. For efficient calculations and a better
understanding of the sign, we then introduce and characterize the sign of single 1-factors. We
show that the product over all signs of all the 1-factors in a 1-factorization is the sign of that
1-factorization. Using this result in an algorithm, we attempt to prove the List Edge Coloring
Conjecture for all graphs with up to 10 vertices. This leaves us with some exceptional cases that
need to be attacked with other methods.
\end{abstract}


\section{Introduction}\label{sec.int}

Using the polynomial method, we prove the \LECC\footnote{See \cite[Section\,12.20]{jeto} for a
discussion of the origins of this coloring conjecture.} for many small graphs $G$. This means, if
such a graph $G$ can be edge colored with $k$ colors ($\chi'(G)\leq k$), then it can also be edge
colored if the color of each edge $e$ has to be taken from an arbitrarily chosen individual list
$L_e$ of $k$ colors ($\chi'_\ell(G)\leq k$). There are no restriction on the lists, apart from the
given cardinality $k$. So, in general, there are very many essentially different list assignments
$e\mto L_e$, and brute-force attempts to find one coloring from every system of lists are
computationally impossible. A way out may be found in the \CN, which seems to be one of our
strongest tools. It can also be used for list coloring of the vertices of a graph (see\,\cite{al}), but it
becomes even more powerful if applied to edge colorings of regular graphs. Ellingham and
Goddyn \cite{elgo} used it to prove the \LECC for regular planar graphs of class\,1\´. As, by
definition, the edges of a class\,1 graph $G$ can be partitioned into $\Delta(G)$ color classes, the
regular class\,1 graphs are precisely the \(1\)"~factorable graphs. \(1\)"~factorable graphs, as we
call regular class\,1 graphs from now on, are also the first target in the current paper,
but our results have implications for other graphs as well. In our previous paper\,\cite{schKp}, 
we could already prove the \LECC for infinitely many \(1\)"~factorable complete graphs. There, we
used a group action in connection with the \CN. Häggkvist and Jansson \cite{haja} could prove the
conjecture for all complete graphs of class\,2\´. Nobody, however, has a proof for $K_{16}$, and
$120$ edges and $15$ colors are completely out of reach for all known numeric methods,
including the algorithms that we suggest here. That we cannot even prove the conjecture for all
complete graphs shows how hard the problem is. Before this background, it is surprising that
Galvin could prove the conjecture for all bipartite graphs \cite{ga}. His proof does not use the \CN,
but the so-called kernel method. Other methods were also used by Kahn \cite{ka}\´, who showed
that the \LECC holds asymptotically, in some sense. Moreover, most of the mentioned results can
also be generalized to edge painting\,\cite{schPC,schPCN}\´, an on-line version of list coloring that
allows alterations of the lists during the coloration process.

This paper has three further sections, and an appendix containing our algorithm. In
Section\,\ref{sec.edp}\´, we formulate a corollary to the \CN that does not require the consideration
of any coefficients of polynomials, only evaluations of polynomial functions. There, we also explain
a well-known connections between the sum of the signs over all 1-factorizations (edge colorings)
of a graph and the List Edge Coloring Conjecture. In Section\,\ref{sec.ItC}, we then provide
another characterization of the sign. We explain how this can be used to calculate the sum of the
signs over all 1-factorizations more efficiently. In Section\,\ref{sec.LCI}\´, we explain to which
conclusions this approach and our computer experiments with graphs on up to $10$ vertices led.


\section{A Nullstellensatz for List Colorings}\label{sec.edp}

We start our investigations from the following coefficient formula \cite{schAlg}\´:

\begin{theorem}[Quantitative Combinatorial Nullstellensatz]\label{sz.cn}\quad

Let $L_1,L_2,\dotsc,L_n$ be finite non"=empty subsets of a field $\F[]$, set
$L:=L_1\times L_2\times\dotsm\times L_n$ and define $d:=(d_1,d_2,\dotsc,d_n)$ via
$d_j:=\abs{L_j}-1$. For polynomials
  $P=\sum_{\d\in\N^n}P_\d x^\d\,\in\,\F[\,][x_1,\dotsc,x_n]$
  of total degree $\deg(P)\leq d_1+d_2+\dotsb+d_n$, we have
\begin{equation*}
  P_d\,=\,\sum_{x\in L}\,N_L(x)^{-1}P(x)\ ,
\end{equation*}
  where $N_L(x)=N_L(x_1,\dotsc,x_n):=\prod_j N_{L_j}(x_j)$ with $N_{L_j}(x_j):=\!\!\!\prod\limits_{\xi\in L_j\!\ssm x_j\!\!\!\!\!}\!\!\!(x_j-\xi)\neq0$\!.

  In particular
,
  if $\deg(P)\leq d_1+d_2+\dotsb+d_n$ then
\begin{equation*}
  P_d\,\neq\,0\ \ \lTo\ \ \ex x\in L\DP P(x)\neq0\ .
\end{equation*}
\end{theorem}

The implication in the second part is known as Alon's Combinatorial Nullstellensatz\,\cite{al2}\´.
The coefficient $P_d$ seems to plays a central role in the Combinatorial Nullstellensatz, but it is
not really important in various applications. One may get a wrong impression form the fact that
$P_d$ is assumed as non-zero in that implication. There are applications of the theorem if the
total degree $\deg(P)$ is strictly smaller than $d_1+d_2+\dotsb+d_n$, and thus $P_d=0$. If
$P_d=0$, then it cannot be that only one summand in the sum in that theorem is non-zero, and
this mens that there cannot be only one solution to the problem that was modeled by $P$\!. So, if
there exist a solution, say a trivial solution, than there must also be a second solution, a non-trivial
solution. This is a very elegant line of reasoning, and it does not require us to look at the
coefficient $P_d$ at all. It is enough to know that the total degree is smaller than
$d_1+d_2+\dotsb+d_n$ and that there is a single trivial solution. Beyond that, the theorem can
also be used to prove the existence of solutions to problems that do not have a trivial solution, for
example the existence of a list coloring of a graph. 
In these cases, looking at the ``leading coefficient'' $P_d$ appears to be unavoidable.
However, to actually calculate $P_d$, usually, the best idea is to use the Quantitative
Combinatorial Nullstellensatz again, just with changed lists $L_j$. In fact, the polynomial $P$ can
be changed, too, as long as the ``leading coefficient'' is not altered. So, theoretically, we can
calculate $P_d$ by applying the theorem to modified lists $\tilde L_j$ and a modified polynomial
$\tilde P$\!. Afterwards, the theorem can then be applied a second time, to $P$ and the original
lists $L_j$, in order to prove the existence of a certain object. In this process, the coefficient
$P_d$ stands in the middle, playing a crucial role. The coefficient $P_d$, however, does not
appear in the initial setting and also not in the final conclusion. Therefore, it must be possible to
formulate a all-in-one ready-to-go corollary in which $P_d$ does not occur. In providing that
corollary, we free the user from the need to understand what $P_d$ is. Of course, in its most
general form, there are two polynomials $P$ and $\tilde P$\!, and two list systems $L$ and $\tilde
L$, which make that corollary look more technical, but it avoids mentioning $P_d$ and should be
easier to apply in many situations:

\begin{corollary}\label{cor.cn}
For $j=1,2,\dotsc,n$, let $L_j$ and $\tilde L_j$ be finite non"=empty subsets of a field
$\F[]$ with $\abs{L_j}=\abs{\tilde L_j}$. 
Let $N_L$ and $N_{\tilde L}$ be the corresponding coefficient functions over the cartesian
products $L$ and $\tilde L$ of these sets. If two polynomials $P,\tilde P\in\F[\,][x_1,\dotsc,x_n]$ of
total degree at most
$\abs{L_1}+\abs{L_2}+\dotsb+\abs{L_n}-n$ have the same homogenous component of degree
$\abs{L_1}+\abs{L_2}+\dotsb+\abs{L_n}-n$
(or at least $\tilde P_d=P_d$), then
\begin{equation*}
  \sum_{x\in\tilde L}\,N_{\tilde L}(x)^{-1}\tilde P(x)\,=\,\sum_{x\in L}\,N_L(x)^{-1}P(x)
\end{equation*}
and, in particular\,\footnote{Also \cite[Th.\,4.5]{schPCN}\´: \ $\sum N_{\tilde L}(x)^{-1}\tilde
P(x)\neq0\ \lTo\textit{$P$\! is \(\,(\abs{\tilde L_1},\dotsc,\abs{\tilde L_n})\)-paintable}$. }\!,
\begin{equation*}
  \sum_{x\in\tilde L}\,N_{\tilde L}(x)^{-1}\tilde P(x)\,\neq\,0\ \ \lTo\ \ \ex x\in L\DP P(x)\neq0\ .
\end{equation*}
\end{corollary}

We want to use this corollary to verify the existence of list colorings of graphs. Therefore, we
apply the corollary to the \emph{edge distance polynomials} $P_{G}$ of graphs $G$. The edge
distance polynomial of a multi-graph $G$ on vertices $v_1,v_2,\dotsc,v_n$ is a polynomial in the
variables $x_1,x_2,\dotsc,x_n,$ with one variable $x_i$ for each vertex $v_i$. It is defined as the
product over all differences $x_i-x_j$ with $v_iv_j\in E(G)$ and $i<j,$ where the factor $x_i-x_j$
occurs as many times in $P$ as the edge $v_iv_j$ occurs in the multi-set $E(G)$. It is also called
the graph polynomial and was introduced in \cite{pe}. We may view it as a polynomial over any
field $\F[].$ If $P_G$ is non-zero at a point $(x_1,x_2,\dotsc,x_n)$ then the assignment $v_i\mto
x_i$ is a proper vertex coloring of $G.$
If the colors $x_i$ are supposed to lie in certain lists $L_i$ then the point $(x_1,x_2,\dotsc,x_m)$
just has to be taken from the Cartesian product $L_1\times L_2\times\dotsm\times L_m.$ Here, we
simple need to assume that the sets $L_i$ lie in $\F[],$ or in an extension field of $\F[].$ This is no
restriction, as one can easily embed the color lists (and their full union $\bigcup_iL_i$) into any big
enough field $\F[].$ We might just take $\F[]=\Q.$ With this ideas our corollary leads to the
following more special result:

\begin{corollary}\label{cor2.cn}
Let $G$ be a multi-graph on the vertices $v_1,v_2,\dotsc,v_n$. To each edge $e$,
between any vertices $v_i$ and $v_j$ with $i<j$, choose a label $a_e$ in a field $\F[]$
(possible $a_e=0$) and associate the monomial $x_i-x_j-a_e$ to the edge $e$. Let
$P$ be the product over all these monomials. For $j=1,2,\dotsc,n$, let $L_j$ be a finite
non"=empty subset of $\F[],$ and define $\ell=(\ell_1,\ell_2,\dotsc,\ell_n)$ via
$\ell_j:=\abs{L_j}$. If $\abs{E(G)}\leq\ell_1+\ell_2+\dotsb+\ell_n-n$ then
\begin{equation*}
  \sum_{x\in L}\,N_L(x)^{-1} P(x)\,\neq\,0
  \ \ \lTo\ \ \text{$G$ is \(\ell\)"~list colorable and \(\ell\)"~paintable.}
\end{equation*}
\end{corollary}

In applications, one will often choose the $a_e$ as zero and take the lists $L_j$ all equal, but
there are also examples where more complicated choices succeeded, as for example in the proof
of the last lemma in \cite{schKp}. Things can be further simplified if we examine edge colorings. In
that case, one has to consider the line graph \rand$"L(G)"$ of $G$ and its edge distance
polynomial $P_{L(G)}$. If $G$ is \(k\)"~regular, then $L(G)$ is the edge disjoint union of $n$
complete graphs $K_k$, and $P_{L(G)}$ factors into $n$ factors accordingly. For each vertex
$v\in V(G)$ there is one complete graph $K_k$ whose vertices are the edges \rand$e\in "E(v_j)"$
incident with $v$. The corresponding factor of $P_{L(G)}$ is the edge distance polynomial
$P_{K_k}(x_e\mit e\in E(v_j))$ of that $K_k$. If the \(k\)"~regular graph is of class\,1, i.e.\ if its
edges can be colored with $k$ colors, then, in the corresponding vertex colorings of $L(G)$, every
color occurs one time at each vertex of that $K_k$. Therefore, by choosing equal lists, say all
equal to \rand$"(k]":=\{1,2,\dotsc,k\}$, the coefficients $N_L(x)^{-1}$ in the sum in the last
corollary become all the same. More precisely, $N_L(x)=N_L(y)$ if $P_{L(G)}(x)\neq0$ and
$P_{L(G)}(y)\neq0$. Moreover, $P_{L(G)}(x)$ assumes, up to the sign, the same value for every
edge coloring $x\DP E(G)\to(k]$. So, in that sum, one basically only has to see which edge
colorings contribute a positive sign and which ones a negative sign. This was already observed in
\cite{al}. It is easy to see that the definition of the sign given there depicts what we need, but we
simplify that a bit. Basically, we only have to be able to say if two edge colorings have same or
opposite sign.
If \randd$"c"\DP E\to(k]$ and \randd$"c_0"\DP E\to(k]$ are proper edge colorings, then
$c|_{E(v)}$ and $c_0|_{E(v)}$ are bijections form the set $E(v)$ of edges at $v\in V(G)$ to $(k]$,
and we set
\begin{equation}
  \sgn_v(c,c_0):=\sgn\bigl(\bigl(c_0|_{E(v)}\bigr)^{-1}\nach\,c|_{E(v)}\bigr)
  \quad\text{and}\quad
  \sgn(c,c_0):=\prod_{v\in V(G)}\sgn_v(c,c_0)\ ,
\end{equation}
where $(c_0|_{E(v)})^{-1}\nach c|_{E(v)}$ is a permutation of $E(v)$ and
$\sgn((c_0|_{E(v)})^{-1}\nach c|_{E(v)})$ is its usual sign. We could have also defined
$\sgn_v(c,c_0)$ as the sign of the inverse permutation
$(c|_{E(v)})^{-1}\nach\,c_0|_{E(v)}$, or as sign of the permutations
$c|_{E(v)}\nach\,(c_0|_{E(v)})^{-1}$ or $c_0|_{E(v)}\nach\,(c|_{E(v)})^{-1}$ in $S_k$.
This is all the same. It is the right definition here, because the sign of a permutation
$\rho$ in $S_k$ is exactly the sign of the edge distance polynomial $P_{K_k}$ of
$K_k$ evaluated at $(\rho_1,\rho_2,\dotsc,\rho_k)$,
\begin{equation}
  \sgn(\rho)=\dfrac{P_{K_k}(\rho_1,\rho_2,\dotsc,\rho_k)}{\,\Abs{P_{K_k}(\rho_1,\rho_2,\dotsc,\rho_k)}\,}\ .
\end{equation}
Hence, we only need to fix one edge coloring $c_0\DP E(G)\to(k]$ and then count how many
colorings $c\DP E(G)\to(k]$ are positive or negative with respect to that \emph{reference
coloring}. It is convenient to define an absolute sign $\sgn(c)$ through
\begin{equation}
  \sgn(c)\,:=\,\sgn(c,c_0)\sgn(c_0)\,,
\end{equation}
where $\sgn(c_0)$ is fixed given as either $+1$ or $-1$. In this section, however, it
does not mater whether $c_0$ is viewed as positive or negative, and we postpone the
stipulation of $\sgn(c_0)$ till later. With that, we arrive at \cite[Corollary\,3.9]{al}:

\begin{corollary}\label{cor3.cn}
Let $G=(V,E)$ be a \(k\)"~regular graph and let \rand$"C(G)"$ be the set of its proper edge
colorings $c\DP E\lto(k]$. Then
\begin{equation*}
  \sum_{c\in C(G)}\,\sgn(c)\,\neq\,0
  \ \ \lTo\ \ \text{$G$ is \(k\)"~list edge colorable and edge \(k\)"~paintable.}
\end{equation*}
\end{corollary}

Actually, we may assume that $G$ has even many vertices, as \(1\)"~factors and \(k\)"~edge
colorings only exist if there are even many vertices. If we exchange two colors in an edge coloring
$c\DP E\to(k]$ of a \(k\)"~regular graph $G$, then all the factors $\sgn_v(c)$ in $\sgn(c)$ change,
but sign $\sgn(c)$ does not change. Therefore, it makes sense to define the sign of a
\emph{\(1\)"~factorization}. A \(1\)"~factorization $F$ of $G$ is a partition
$F=\{F_1,F_2,\dotsc,F_k\}$ of the edge set $E(G)$ into $k$ \(1\)"~factors (perfect matchings). To
every \(1\)"~factorization $F$ there are $k!$ edge colorings $c$ with $F$ as set of fibers
$c^{-1}(\{\ä\})$. All of them have the same sign, and we define
\begin{equation}
\sgn(F):=\sgn(c)\ .
\end{equation}
With that, the last corollary can be rewritten as follows:

\begin{corollary}\label{cor4.cn}
Let $G=(V,E)$ be a \(k\)"~regular graph 
and let \rand$"\OF(G)"$ be the set of \(1\)"~factorizations of $G$. Then
\begin{equation*}
  \sum_{F\in\OF(G)}\,\sgn(F)\,\neq\,0
  \ \ \lTo\ \ \text{$G$ is \(k\)"~list edge colorable and edge \(k\)"~paintable.}
\end{equation*}
\end{corollary}



\section{Another Characterization of the Sign}\label{sec.ItC}

In this section, $G$ denotes a \(k\)"~regular graph on the vertices
$v_1,v_2,\dotsc,v_{2n}$, and $F=\{F_1,F_2,\dotsc,F_k\}$ denotes a \(1\)"~factorization
of $G$. We examine the sign $\sgn(F)$ in more detail, starting from the following
definition:

\begin{definition}
Let $F_1=\{e_1,e_2,\dotsc,e_n\}$ be a \(1\)"~factor of a \(k\)"~regular graph $G$ on the vertices
$v_1,v_2,\dotsc,v_{2n}$. Let $1\leq i_k<j_k\leq 2n$ be such that $e_k=v_{i_k}v_{j_k}$, for
$k=1,2,\dotsc,n$. We say that an edge $e_k\in F_1$ intersects another edge $e_\ell\in F_1$ if
$i_k<i_\ell<j_k<j_\ell$ or $i_\ell<i_k<j_\ell<j_k$. We define
$$
\intt(e_k,e_\ell)\,:=\,\begin{cases}
1 & \text{if $e_k$ intersects $e_\ell$,}\\
0 & \text{otherwise,}
\end{cases}
$$
and set
$$
\intt(F_1)\,:=\,\sum_{1\leq k<\ell\leq n}\intt(e_k,e_\ell)\quad\ \text{and}\quad
\sgn(F_1)\,:=\,(-1)^{\intt(F_1)}\,.
$$
\end{definition}

If we position the $2n$ vertices consecutively around a cycle and draw the edges as
strait lines, then an intersection is an actual intersection between lines. With this picture
in mind, it is not hard to see that, if $\intt(v_iv_j,F_1)$ denotes the number of
intersections of an edge $v_iv_j\in F_1$ with other edges in $F_1$, then
\begin{equation}\label{eq.cuts}
  \intt(v_iv_j,F_1)\,\equiv\,j-i-1\pmod{2}\ .
\end{equation}
This, however, does not help to determine the sign $\sgn(F_1)$ of $F_1$, as
\begin{equation}
  \sum_{e\in F_1}\intt(e,F_1)\,=\,2\intt(F_1)\ ,
\end{equation}
with a $2$ in front of $\intt(F_1)$. Counting the number of all intersections of each edge $e$ is not
the right approach here. We may order $F_1$ to $\overrightarrow{F_1}=(e_1,e_2,\dotsc,e_n)$
and count only the intersections of an edge $e_k$ with the \emph{subsequent edges} $e_\ell$, $\ell>k$. 
If $\intt(e_k,\overrightarrow{F_1})$ denotes this number, then the corresponding sum yields the
desired result,
\begin{equation}
  \intt(F_1)\,=\,\sum_{e\in F_1}\intt(e,\overrightarrow{F_1})\ .
\end{equation}
Hence,
\begin{equation}
  \sgn(F_1)\,=\,\prod_{e\in F_1}\sgn(e,\overrightarrow{F_1})\ ,
\end{equation}
if we set
\begin{equation}
  \sgn(e,\overrightarrow{F_1})\,:=\,(-1)^{\intt(e,\overrightarrow{F_1})}\ .
\end{equation}
This formula may be used to calculate the sign of a \(1\)"~factor in algorithms that generate a
\(1\)"~factor by successively adding single edges. And, there is also an analog to
Formula\,\eqref{eq.cuts}\´. We may just count how many of the vertices $b$ that lie between the
two ends $v_{i_k}$ and $v_{j_k}$ of the edge $e_k$ are not yet matched when we add $e_k$ to
the sequence $(e_1,e_2,\dotsc,e_{k-1})$. So,
\begin{equation}\label{eq.cuts2}
  \intt(v_{i_k}v_{j_k},\overrightarrow{F_1})\,\equiv\,\Abs{\{\´b\!\mit\!i_k<b<j_k\,,\,\,b\nin e_1\cup e_2\cup\dotsb\cup e_{k-1}\}}\pmod{2}\ .
\end{equation}
In our algorithm, we kept track of these unmatched $b$ by using a doubly linked linear lists. From
each unmatched vertex $b$, we have at any time a link to the unmatched vertex before $b$ and a
link to the unmatched vertex after $b$. Updating these links can then be done without shifting all
subsequent vertices one place forward.

The next theorem shows that the signs of the \(1\)"~factors in a \(1\)"~factorization $F$ can be
used to calculate the sign of $F$\!. This can then be used in algorithms that calculate the
\(1\)"~factorizations of a graph by successively adding new \(1\)"~factors. The advantage is that
the sign of a \(1\)"~factor that is added at a certain point has to be calculated only once, for all the
\(1\)"~factorizations that are generate afterwards, by adding more \(1\)"~factors in all possible
ways. It is clear that the formula in the next theorem 
does not really depend on the sign of the underlying reference coloring $c_0$, or the equivalent
\emph{reference \(1\)"~factorization} $\bigl\{c_0^{-1}(\{\ä\})\mit\ä\in(k]\bigr\}$. But, to avoid
additional minus signs in the theorem, we synchronize our different signs at this point, and define
 \rand\begin{equation}
  "\sgn(c_0)"\,:=\,\prod_{\ä\in(k]}(-1)^{\intt(c_0^{-1}(\{\ä\}))}\,=\,(-1)^{\intt(c_0)}\,\in\,\{-1,+1\}\,,
\end{equation}
where
 \rand\begin{equation}
"\intt(c_0)":=\sum_{\ä\in(k]}\intt(c_0^{-1}(\{\ä\}))
\end{equation}
is the number of intersections between edges of equal color in $c_0$, if the vertices
$v_1,v_2,\dotsc,v_{2n}$ are arranged consecutively on a cycle and the edges are drown as strait
lines. With this stipulation of the sign of the reference coloring $c_0$, we have the following
theorem:


\begin{theorem}
Let $G=(V,E)$ be a \(k\)"~regular graph on the vertices $v_1,v_2,\dotsc,v_{2n}$, and
let $F=\{F_1,F_2,\dotsc,F_k\}$ be a \(1\)"~factorization of $G$. Then
\begin{equation*}
  \sgn(F)\,=\,\prod_{i=1}^k\sgn(F_i)\ .
\end{equation*}
In other words, if $c\DP E\lto (k]$ is an edge coloring, then
\begin{equation*}
  \sgn(c)\,=\,(-1)^{\intt(c)}\ ,
\end{equation*}
where $\intt(c)$ is the number of intersections between edges of equal color, if the vertices
$v_1,v_2,\dotsc,v_{2n}$ of $G$ are arranged consecutively on a cycle and the edges are drown
as strait lines.
\end{theorem}

Proving this theorem is the main task of this section. We do this in a tropologic way, using
Jordan's Curve Theorem. From this theorem, we know that any two closed curves on the sphere
have even many intersections with each other (and that even if they also have intersection points
with them selves, which we just do not count). We also use the fact that the sign of a permutation
$\rho\in S_k$ is $-1$ to the power of the number of \emph{inversions} of $\rho$. Here, a pair
$(i_1,i_2)\in(k]^2$ with $i_1<i_2$ is an inversion of $\rho$ if $\rho(i_1)>\rho(i_2)$. We will use that
this property can be characterized as intersection of strait lines in $\Rl^2$\!. Indeed, the pair
$(i_1,i_2)$ is an inversion if and only if the line from $(i_1,h_1)$ to $(\rho(i_1),h_2)$ intersects
with the line from $(i_2,h_1)$ to $(\rho(i_2),h_2)$, where $h_1$ and $h_2$ are any two different
real numbers:

\begin{proof}
Let $c_0\DP E\to(k]$ be the reference coloring of $G$, and let $c\DP E\to(k]$ be another edge
coloring. We have to show that $\intt(c)\equiv\intt(c_0)\pmod{2}$ if and only if $\sgn(c,c_0)=1$. To
compare the numbers of intersections in $c$ and $c_0$, we draw the colored graph $(G,c)$ on
top of a round cylinder, with the vertices in counter-clockwise order along the boundery of the
upper disc. The colored graph $(G,c_0)$ is drawn on the bottom of the cylinder, in such a way
that every vertex $v_j$ of $(G,c)$ lies vertically above the corresponding vertex $v_j$ of
$(G,c_0)$. Now, we remove the vertex $v_j$ in $(G,c)$ and $(G,c_0)$ and connect the open ends
of the edges in $E(v_j)$ on the top disk with those in the bottom disk. We connect edges of equal
color by a line along the lateral surface of the cylinder. As to every color $\alpha\in(k]$ there exists
exactly one edge of color $\alpha$ incident with $v_j$ in $(G,c)$ and in $(G,c_0)$, this makes
exactly one line of every color (for every $j\in(2n]$). To avoid that these $k$ lines lie on top of
each other, we assume that we have cut down the radius of the cylinder a bit, so that the edges in
$E(v_j)$ do not end in exactly the same point of the boundary of the upper, resp.\ lower, disc.
Hence, we have $2n$ disjoint intervals $I_j$ on the edge of each disc, corresponding to the $2n$
removed vertices $v_j$. In each interval $I_j$, on each disc, the edges of $E(v_j)$ arrive in
consecutive order, corresponding to the clockwise order of the edges in $E(v_j)$ around $v_j$.
We may imagine the area between the upper interval $I_j$ and the lower interval $I_j$ as a
rectangle with $k$ straight but slanted lines crossing from the upper interval $I_j$ to the lower
interval $I_j$. If a color $\ä$ occurs on, say, the $2^{\text{nd}}$ edge of $E(v_j)$ in $c$, and on,
say, the $5^{\text{th}}$ edge of $E(v_j)$ in $c_0$, then there is a line of color $\ä$ running from
the $2^{\text{nd}}$ position in the upper interval $I_j$ to the $5^{\text{th}}$ position in the lower
interval $I_j$.\smallskip

\emph{Claim:} 
$\sgn_{v_j}(c,c_0)$ is equal to $-1$ to the power of the
number of intersections between the $k$ lines that run from the upper interval $I_j$ to
the lower interval $I_j$. 
$\sgn(c,c_0)$ is equal to $-1$ to the power of
the number of intersections between all lines on the lateral surface of the
cylinder.\smallskip

We prove the first part of this claim by observing that every intersection corresponds to
an inversion of the permutation $\rho:=\bigl(c_0|_{E(v_j)}\bigr)^{-1}\nach\,c|_{E(v_j)}$ of
$E(v_j)$. We identify the clockwise ordered edges $e_1,e_2,\dotsc,e_k$ in $E(v_j)$,
and the position in $I_j$ where they arrive, with the integers $1,2,\dotsc,k$ (in that
order). With that identification, $\rho$ is actually an element of $S_k$, and the pair
$(1,2)$, for instance, is an inversion of $\rho$ if and only if the lateral lines that start in
position $1$ and $2$ of the upper interval $I_j$ intersect. Obviously, the colors of these
two lines are $c(e_1)$ and $c(e_2)$, respectively. Inside $(G,c_0)$, these two colors
occur at the edges $\rho(e_1)=\bigl(c_0|_{E(v_j)}\bigr)^{-1}(c(e_1))$ and $\rho(e_2)$ of
$E(v_j)$, respectively. So, position $1$ and $2$ in the upper interval $I_j$ are
connected to position $\rho(e_1)$ and $\rho(e_2)$ in the lower interval $I_j$. Our two
lines cross if and only if $\rho(e_1)>\rho(e_2)$, if and only if $(1,2)$ is an inversion of
$\rho$. The first part of our claim follows from that. It holds for each $j\in(2n]$, and that
is just summed up in the second part.\smallskip

From the claim, we see that $\sgn(c,c_0)=1$ if and only if the number of intersections between
lines on the lateral surface of the cylinder is even. Note also that all these lateral intersections are
intersections between lines of different color. Overall, on the whole cylinder, there are even many
intersections between lines of different color. This follows from Jordan's Curve Theorem, as all
lines together form a system of monochromatically colored closed curves on the surface of the
cylinder. Therefore, modulo $2$, the number of intersections of differently colored edges in the
upper disk is equal to that number in the lower disk if and only if $\sgn(c,c_0)=1$. Since the total
number of intersections (that between differently and equally colored edges) is the same on both
disks, this also means that $\intt(c)\equiv\intt(c_0)\pmod{2}$ if and only if $\sgn(c,c_0)=1$.
\end{proof}


\section{The List Chromatic Index of Small Graphs}\label{sec.LCI}

Based on Corollary\,\ref{cor4.cn} and the results of the previous section, we have tried to
determine the list chromatic index $\chi'_{\ell}(G)$ of all graphs on up to 10 vertices, in an attempt
to prove the \LECC for small graphs. We implemented the approach explained in the previous
sections in SageMath \cite{sage}\´, importing regular graphs from the webpage \cite{me}
described in \cite{me2}. With that we attacked all regular graphs on 4, 6, 8 or 10 vertices. The
results are shown in the first paragraph of the following subsection. We tried than to draw
conclusions about the list chromatic index of all graphs with up to 10 vertices. We did this by
considering embeddings into regular graphs on even many vertices. Unfortunately, there are
many exceptional cases and special circumstances. We report about these difficulties, and some
ideas how to overcome them, in quite a view case distinctions. It was not possible to go through
all the cases and to prove the \LECC for all graphs on up to 10 vertices. If, however, someone
wants to prove the \LECC for just one particular small graph, he or she may find a way to do so
within our case distinctions.

In the following case distinctions, the word \emph{graph} stands for connected graph, and a
regular graph $G$ is a \emph{zero-sum graph} if the sum $\sum\sgn(F)$ over all
\(1\)"~factorizations $F\in\OF(G)$ vanishes. We call a graph \emph{small} if it has at most $10$
vertices, and we call it \emph{even} resp.\ \emph{odd} if it has even resp.\ odd many vertices.


\subsection{Small even graphs}\label{sec.seg}

\paragraph{Regular Graphs.}
By checking all small regular even graphs, we found only three graphs of class\,2\´. The Petersen
graph and the following two graphs:

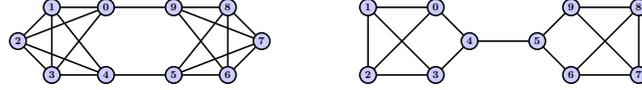
\begin{center}
\begin{tikzpicture}
  [scale=.9,line width=0.6pt,every node/.style={draw,circle,fill=blue!20,font=\bf,scale=0.4}]
  \node (0) at (1.3,1) {\!0\!};
  \node (1) at (0.5,1) {\!1\!};
  \node (2) at (0,0.5)  {\!2\!};
  \node (3) at (0.5,0)  {\!3\!};
  \node (4) at (1.3,0) {\!4\!};
  \node (5) at (2.3,0) {\!5\!};
  \node (6) at (3.1,0) {\!6\!};
  \node (7) at (3.6,0.5) {\!7\!};
  \node (8) at (3.1,1) {\!8\!};
  \node (9) at (2.3,1) {\!9\!};
    \path[every node/.style={font=\sffamily\small,scale=0.8,font=\bf,fill=white,sloped,thick}]
      (0) edge (1)
      (0) edge (2)
      (0) edge (3)
      (0) edge (9)
      (1) edge (2)
      (1) edge (3)
      (1) edge (4)
      (2) edge (3)
      (2) edge (4)
      (3) edge (4)
      (4) edge (5)
      (5) edge (6)
      (5) edge (7)
      (5) edge (8)
      (6) edge (7)
      (6) edge (8)
      (6) edge (9)
      (7) edge (8)
      (7) edge (9)
      (8) edge (9);
\end{tikzpicture}
\qquad\quad
\begin{tikzpicture}
  [scale=.9,line width=0.6pt,every node/.style={draw,circle,fill=blue!20,font=\bf,scale=0.4}]
  \node (0) at (1.5,1) {\!0\!};
  \node (1) at (0.5,1) {\!1\!};
  \node (2) at (2,0.5)  {\!4\!};
  \node (3) at (0.5,0)  {\!2\!};
  \node (4) at (1.5,0) {\!3\!};
  \node (5) at (3.5,0) {\!6\!};
  \node (6) at (4.5,0) {\!7\!};
  \node (7) at (3,0.5) {\!5\!};
  \node (8) at (4.5,1) {\!8\!};
  \node (9) at (3.5,1) {\!9\!};
    \path[every node/.style={font=\sffamily\small,scale=0.8,font=\bf,fill=white,sloped,thick}]
      (0) edge (1)
      (0) edge (2)
      (0) edge (3)
      (1) edge (3)
      (1) edge (4)
      (2) edge (4)
      (2) edge (7)
      (3) edge (4)
      (5) edge (6)
      (5) edge (7)
      (5) edge (8)
      (6) edge (8)
      (6) edge (9)
      (7) edge (9)
      (8) edge (9);
\end{tikzpicture}
\end{center}
Our main method does not apply to class\,2 graphs. In these three cases, however, one can
simply add a suitable \(1\)"~factor, and prove the \LECC for the resulting graph of class\,1\´. It is,
in fact, possible to choose the \(1\)"~factor in a way that the extended graph is not a zero-sum
graph. So, in the shown three cases, the \LECC holds. Unfortunately, our method also failed in a
number of other cases, where the sum $\sum\sgn(F)$ over all \(1\)"~factorizations $F\in\OF(G)$
simple was zero. The smallest zero-sum graph is $K_{3,3}$, but this graph is bipartite. Hence, it
meets the \LECC by Galvin's Theorem \cite{ga}\´. On $8$ vertices, there are exactly three
zero-sum graphs. The complement $\overline{C_3\cup C_5}$ of the disjoint union of a \(3\)"~cycle
and a \(5\)"~cycle, and the following graph and its complement:
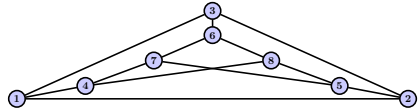
\begin{center}
\begin{tikzpicture}
  [rotate=0,y=4mm,line width=0.6pt,scale=1.3,every node/.style={draw,circle,fill=blue!20,font=\bf,scale=0.4}]
  \node (0) at (2,1.63)  {\!6\!};
  \node (1) at (1.4,0.98)  {\!7\!};
  \node (2) at (4,0) {\!2\!};
  \node (3) at (3.3,0.33) {\!5\!};
  \node (4) at (0.7,0.33) {\!4\!};
  \node (5) at (0,0) {\!1\!};
  \node (6) at (2.6,0.98) {\!8\!};
  \node (7) at (2,2.25) {\!3\!};
    \path[every node/.style={font=\sffamily\small,scale=0.8,font=\bf,fill=white,sloped}]
      (0) edge (1)
      (2) edge (7)
      (3) edge (6)
      (4) edge (5) %
      (0) edge (6)
      (1) edge (4)
      (2) edge (3)
      (5) edge (7) %
      (0) edge (7)
      (1) edge (3)
      (2) edge (5)
      (4) edge (6); 
\end{tikzpicture}
\end{center}
On $10$ vertices there are $51$ zero-sum graphs out of $164$ regular class\,1 graphs
(\(1\)"~factorable graphs). There are $5$ zero-sum graphs of degree $3$, $17$ of degree $4$,
$18$ of degree $5$, $8$ of degree $6$, and $3$ of degree $7$. It seems that, in every small
zero-sum graph, one can find a symmetry of order $2$ that turns even edge coloring ($\sgn=+1$)
into odd ones ($\sgn=-1$) and vice versa; which explains the vanishing sum. The most simple
symmetry of this kind is given if two non-adjacent vertices of odd degree have the same
neighbors, or if two adjacent vertices of even degree have the same neighbors. But, there are
also more complicated cases. In the complement of the Petersen graph, for example, it is more
difficult to understand how odd and even edge colorings are matched through a graph symmetry.
Overall, it should be possible to proof the \LECC for all found zero-sum graphs with other
methods. Some well chosen case distinctions with respect to the color lists might suffice. This kind
of reasoning, however, is usually quite tedious and depends very much on the structure of the
graph.

\paragraph{Non-regular Graphs.}
If a regular graph $G$ is of class\,1 and meets the \LECC, then every subgraph of same maximal
degree still is of class\,1 and still meets the \LECC. With this argument, most non-regular small
even graphs can be proven to be of class\,1 and to meet the \LECC. We just have to consider
regular even extensions of same maximal degree. If an extension is still small, we may apply our
findings about small regular even graphs. There are, however, three difficulties:
\smallskip\\
(i) Some small non-regular even graphs cannot be embedded into a regular graph by adding
edges only, which would keep these graphs small. Several examples of this kind can be
constructed from \(k\)"~regular graphs ($k\geq3$) that contain an induced path $u{-}v{-}w$ by
removing the edges $uv$ and $vw$\!, and inserting the edge $uw$\!.
\smallskip\\
(ii) The three small regular even graphs of class\,2 are not suitable as regular extensions in this
line of reasoning. Some of their subgraphs are actually of class\,2\´,
and we can only conclude that these class\,2 subgraphs meet the \LECC.
\smallskip\\
(iii) There are still some open cases among the small regular even class\,1 graphs, for which we
not yet have proven the conjecture. Circumventing these cases is not always possible, as there
may not be many different ways to add edges.

\subsection{Small odd graphs}\label{sec.sog}

\paragraph{Class\,2 Graphs (including all Regular Graphs).}
All regular graphs of odd order are of class\,2\´, as no \(1\)"~factors exist. Moreover, if we start
from an \(k\)"~regular odd graph and remove less than $k/2$ edges, then the graph remains in
class\,2\´, because it is still \emph{overfull} ($\abs{E}>\Delta\cdot\cop{\abs{V}/2}$). All graphs that
we obtain in this way have maximal degree $k$, which is necessarily an even number, as the
initial regular graph was odd. Odd class\,2 graphs with odd maximal degree are not obtained in
this way. But, they do exist. One example is $K_8$ with one edge subdivided by a new vertex,
which is still overfull.
To prove the \LECC for this graph and for all class\,2 graphs $G$\!, however, we do not need to
embed $G$ into a regular class\,2 graph of same maximal degree $\Delta(G)$\!. To prove that a
graph $G$ (whether of class\,2 or not) has list chromatic index $\Delta(G)+1$, we may simply
embed it into a class\,1 graph whose maximal degree is $\Delta(G)+1$. If the \LECC was proven
for that extension graph, then $\chi'_{\ell}(G)\leq\Delta(G)+1$, and then the \LECC holds for $G$ if
$G$ is of class\,2\´. We may also add vertices. In this way, most small odd graphs can
be embedded into a suitable regular graph. 
As in the case of even non-regular graphs, however, there are three difficulties:
\smallskip\\
(i) Some small odd graphs cannot be embedded into a regular graph by adding only one vertex
and some edges, which would keep these graphs small. One example of this kind is $K_8$ with
one edge subdivided by a new vertex.
\smallskip\\
(ii) The three small regular even graphs of class\,2 are not suitable as regular extensions in this
line of reasoning and must be circumvented. Since the maximal degree can go up by one,
however, there is a lot of flexibility. One can show that the three exceptions of class\,2 are not
needed as extension graphs. Still, circumventing them is an additional difficulty if one tries to draw
general conclusions.
\smallskip\\
(iii) There are still some open cases among the small regular even class\,1 graphs. If we try to
embed a single small odd class\,2 graph, it is often easy to circumvent the open cases. But, in
general examinations, avoiding open cases is difficult.

\paragraph{Class\,1 graphs.}
The majority of small odd graphs are of class\,1 and, in particular, non-regular. For these graphs,
embedding without increasing the maximal degree frequently works. One can try to add just one
vertex and some additional edges. In this way, the results about small even regular graphs can be
applied. As in the other case where we discussed embedding, 
there are three difficulties:
\smallskip\\
(i) Adding just one vertex, to stay within the small graphs, does not work if there are not enough
vertices of sub-maximal degree to which the new vertex can be connected. In this regard, there
are obviously more problematic cases as in the discussion of small odd non-regular graphs of
class\,2\´, where we could increase the maximal degree by one.
\smallskip\\
(ii) The three small regular even graphs of class\,2 are not suitable as regular extensions in this
line of reasoning. However, if we remove just one vertex from any of them, they remain in
class\,2\´. Hence, the three class\,2 graphs do not appear as single-vertex extensions of class\,1
graphs. And, if we need to add a vertex plus some edges, we may be able to circumvent these
three graphs.
\smallskip\\
(iii) If we try to embed a single small odd class\,1 graph, circumventing the open cases among the
small regular even class\,1 graphs is sometimes not possible.

\section{Appendix}

We implemented our algorithm in SageMath \cite{sage} as function $\verb+weighted_sum()+\!$,
useing only commands available in the underlying programming language python.
Equation\,\eqref{eq.cuts2} in Section\,\ref{sec.ItC} provides the foundation for the accumulation of
the sign
\begin{equation}
\sgn(e_k,\overrightarrow{F_1})\,:=\,(-1)^{\intt(e,\overrightarrow{F_1})}
\end{equation}
of an edge $e_k$ with respect to a partial \(1\)"~factor
$\overrightarrow{F_1}=(e_1,e_2,\dotsc,e_{k-1})$ in the variable $\verb+sgn+$. The list
$\verb+previous_Unmatched+$ resp.\ $\verb+next_Unmatched+$ in $\verb+weighted_sum()+$
contains in the cell with number $\verb+um+$ the link to the unmatched vertex before
$\verb+um+$ resp.\ after $\verb+um+$, as explained after Equation\,\eqref{eq.cuts2}\´. By default
these lists are set to $\verb+[-1..9]+$ resp.\ $\verb+[1..11]+$, but they can also be entered as
optional parameters of $\verb+weighted_sum()+\!$. The last entry $\verb+next_Unmatched[-1]+$
of $\verb+next_Unmatched+$ usually points to the very first unmatched vertex. If it is greater or
equal to the number of vertices, however, it means that a fresh bootstrapping needs to be
initiated. This is done in the $\verb+elif+$ part of the initial bootstrapping mechanism in our
recursive function. Here, we also force the first edge of vertex $0$ to be in the first 1-factor, its
second edge to be in the second 1-factor, etc. Hence, in a \(k\)"~regular graph, from the $k!$
equivalent edge colorings that arise out of one edge coloring by permutation of colors, only one is
counted. In out algorithm, we also do not take first the product over all signs of all edges in an
\(1\)"~factorization, to afterward add up all the products that we get for the different
\(1\)"~factorizations. Instead, based on the distributive law, we take the sum over partial
\(1\)"~factorizations during the construction process, and then multiply these partial sums with the
sign of the edge that extends all these partial \(1\)"~factorizations. This speeds up our algorithm.

As input a regular graph on even many vertices is required. The format has to be as in the
example of $K_6$, shown below the definition of $\verb+weighted_sum()+$ in line\,35\´. The lists
of adjacencies of each vertex has to be in strictly increasing order, without showing predeceasing
vertices. For graphs with more than $10$ vertices, the lists $\verb+[-1..9]+$ and $\verb+[1..11]+$
in line 1, 2, 12 and 13 have to be extended. As output, we obtain the sum $\sum\sgn(F)$ over all
\(1\)"~factorizations $F$ of the graph, as it is needed in Corollary\,\ref{cor4.cn}\´:\bigskip

\noindent{\bf Algorithm}\smallskip
\begin{lstlisting}[basewidth={0.58em,0.45em},fontadjust,numbers=left,numberstyle=\tiny,numbersep=8pt]
def weighted_sum(Graph, previous_Unmatched = [-1..9], \
                 next_Unmatched = [1..11]): # 2 optional param.
    # by default, start = next_Unmatched[-1] = 11 > len(Graph)
    # next_Unmatched[j] is the unmatched vertex after j
    # previous_Unmatched[j] is the unmatched vertex before j
    to_match = next_Unmatched[-1] # next_Unmatched[-1] is start
    if to_match < len(Graph): # 1-factor under construction
        neighbors = Graph[to_match]
    elif len(Graph[0]) <> 0: # start next 1-factor
        to_match = 0 # 0 shall be matched first
        neighbors = [Graph[0][0]] # to avoid color permutations
        previous_Unmatched = [-1..9] # fresh bootstrapping
        next_Unmatched = [1..11]
    else: return 1 # 1-factorization complete, edgeless graph
    um = next_Unmatched[to_match]
    previous_Unmatched[um] = -1 # bypass to_match
    next_Unmatched[-1] = um # bypass to_match
    w_sum = 0 # subtotal of weighted_sum()
    sgn = 1 # initial sign of edge {to_match,nbr}
    for i in range(len(neighbors)):
        nbr = neighbors[i] # i^th neighbor of to_match
        while um < nbr: # um is bridged by {to_match,nbr}
            sgn = -sgn # bridged unmatched vertices flip sgn
            um = next_Unmatched[um]
        if um == nbr: # match to_match with nbr
            gr = [[n for n in lst] for lst in Graph] # deepcopy
            del gr[to_match][i] # remove edge {to_match,nbr}
            p_um = [n for n in previous_Unmatched] # deepcopy
            n_um = [n for n in next_Unmatched] # deepcopy
            p_um[n_um[nbr]] = p_um[nbr] # bypass nbr
            n_um[p_um[nbr]] = n_um[nbr] # bypass nbr
            w_sum = w_sum + sgn * weighted_sum(gr,p_um,n_um)
    return w_sum # output w_sum

graph = [[1,2,3,4,5],[2,3,4,5],[3,4,5],[4,5],[5],[]] # K6
# 0 adjacent to 1,2,3,4,5; 1 adjacent to 2,3,4,5 (and 0); etc.
weighted_sum(graph) # the initial call of weighted_sum()
# returns the sum of all signs of all 1-factorizations of graph
\end{lstlisting}


\begin{thebibliography}{6}
\bibitem
{al}
  N.\,Alon: \textit{Restricted Colorings of Graphs.}\\
  In: Surveys in combinatorics, 1993.
  London Math.\ Soc.\ Lecture Notes Ser.\ 187,\\
  Cambridge Univ.\ Press, Cambridge 1993, 1-33.
\bibitem
{al2}
  N.\,Alon: \textit{Combinatorial Nullstellensatz.}\\
  Combin.\ Probab.\ Comput.\ 8, No.\ 1-2 (1999), 7-29.
\bibitem
{elgo}
  M.\,N.\,Ellingham, L.\,Goddyn: \textit{List Edge Colourings of Some 1-Factorable Multigraphs.}   Combinatorica 16 (1996), 343-352.
\bibitem
{ga}
    F.\,Galvin: \textit{The List Chromatic Index of a Bipartite Multigraph.}\\
    J.\,Combin.\ Theory Ser.\ B 63 (1995), 153-158.
\bibitem
{haja}
    R.\,H\"{a}ggkvist  and J.\,Janssen: \textit{New Bounds on the List-Chromatic Index
    of the Complete Graph and Other Simple Graphs.}\\
    Combin.\ Probab.\ Comput.\ 6 (1997), 295-313.
\bibitem
{jeto}
  T.\,R.\,Jensen, B.\,Toft: \textit{Graph Coloring Problems.}   Wiley, New York 1995.
\bibitem
{ka}
  J.\,Kahn: \textit{Asymptotically Good List-Colorings.}\\
  J.\ Comb.\ Theory, Ser.\ A 73(1) (1996), 1-59.
\bibitem
{me}
  M.\,Meringer: \textit{Connected Regular Graphs.}\\
  http://www.mathe2.uni-bayreuth.de/markus/reggraphs.html
\bibitem
{me2}
  M.\,Meringer: \textit{Fast Generation of Regular Graphs and Construction of Cages.}
  Journal of Graph Theory 30 (1999), 137-146.
\bibitem
{pe}
  J.\,Petersen: \textit{Die Theorie der regularen Graphs.}   Acta Math.\ 15 (1891), 193-220.
\bibitem
{sage}
  SageMath, the Sage Mathematics Software System (Version 7.4.1),
  The Sage Developers, 2017, http://www.sagemath.org.
\bibitem
{schAlg}
  U.\,Schauz: \textit{Algebraically Solvable Problems:\\
  Describing Polynomials as Equivalent to Explicit Solutions.}\\
  The Electronic Journal of Combinatorics 15 (2008), \#R10.
\bibitem
{schPC}
  U.\,Schauz: \textit{Mr.\ Paint and Mrs.\ Correct.}\\
  The Electronic Journal of Combinatorics 15 (2008), \#R145.
\bibitem
{schPCN}
  U.\,Schauz: \textit{A Paintability Version of the Combinatorial Nullstellensatz,\\
  and List Colorings of \(k\)-partite \(k\)-uniform Hypergraphs.}\\
  The Electronic Journal of Combinatorics 17/1 (2010), \#R176.
\bibitem
{schKp}
  U.\,Schauz: \textit{Proof of the List Edge Coloring Conjecture for Complete Graphs of Prime Degree.}
  The Electronic Journal of Combinatorics 21/3 (2014), \#P3.43.
\end{thebibliography}
\end{document}